\documentclass[12pt]{amsart}
\usepackage{amsthm}
\usepackage{graphicx}
\usepackage{}
\usepackage{amsfonts}
\usepackage{amssymb}

\newtheorem{theorem}{Theorem}

\newtheorem*{conclusion}{Conclusion}

\newtheorem{corollary}{Corollary}

\newtheorem{definition}{Definition}
\newtheorem{example}{Example}

\newtheorem{lemma}{Lemma}

\newtheorem{proposition}{Proposition}
\newtheorem{remark}{Remark}

\begin{document}

\title[On the cohomology of frobeniusian model Lie algebras]{On the cohomology of frobeniusian model Lie algebras}

\author[J. M. Ancochea and R. Campoamor]{J M Ancochea-Berm\'{u}dez\dag, R Campoamor\dag}

\address{\dag \ Departamento de Geometr\'{\i}a y Topolog\'{\i}a, Facultad CC. Matem\'{a}ticas U.C.M, Avda. Complutense s/n, E-28040 Madrid ( Spain ). \newline Jose\_Ancochea@mat.ucm.es and rutwig@nfssrv.mat.ucm.es}

\begin{abstract}
We compute the first and second cohomology groups with coefficients in the adjoint module of frobeniusian model algebras whose parameters move in a dense open subset of $\mathbb{C}^{p-1}$, and obtain upper bounds for the dimension of cohomology groups of frobeniusian Lie algebras. Moreover, it is shown that for a dense open subset of $\mathbb{C}^{p-1}$ the deformations of model algebras also belong to the family. Therefore any frobeniusian non-model algebra contracts on some element of the model whose parameters move on a finite union of hyperplanes. Further applications as the nullity of Rim's  quadratic map $sq_{1}$ are obtained.
\end{abstract}

\maketitle

\section{Introduction}

Lie algebra model theory arises from the attempts to describe certain neighbourhoods of Lie algebras in the variety of Lie algebra laws $\mathcal{L}^{n}$ in terms of particular properties. While such properties as solvability are nonstable for deformations, other as semisimplicity are preserved, to such an extent that it implies rigidity. There are intermediate properties, usually extracted from Differential Geometry, and which play an important role in the study of forms and differential systems, as for example the existence of linear contact or symplectic forms or $r$-contact systems \cite{G3}. Here we will consider the property of having a linear form whose differential is symplectic. These algebras, which are solvable and are usually called frobeniusian, were classified up to contraction by Goze in \cite{G2}. Now the model is not unique for the described algebras: there exists a family on $\left(p-1\right)$ parameters, which shows that models need not to be isolated. Most of these objects have been studied in the frame of perturbation theory \cite{G1}. This approach allows some simplifications in cohomological questions, such as the deformation equation of Nijenhuis and Richardson \cite{NR}, which reduces to analyze a finite system.\newline From the physical point of view, solvable Lie algebras often occur as Lie algebras of symmetry groups of differential equations, which, added to the symplectic structure, of extraordinary importance, justifies the interest of analyzing frobeniusian Lie algebras.\newline In this paper we calculate the first and second cohomology groups ( with coefficients in the adjoint module) of frobeniusian model algebras whose parameters move in a dense open subset of $\mathbb{C}^{p-1}$. This allows to obtain upper bounds for the dimension of cohomology groups of frobeniusian Lie algebras. Moreover, it is shown that for a dense open subset of $\mathbb{C}^{p-1}$ the deformations of model algebras whose parameters belong to the subset also belong to the family. Therefore any frobeniusian non-model algebra contracts on some element of the model whose parameters move on a finite union of hyperplanes of $\mathbb{C}^{p-1}$. Further applications as the nullity of Rim's  quadratic map $sq_{1}$ are obtained.

\bigskip

Any Lie algebra considered in this work is finite dimensional over the field $\mathbb{C}$. Moreover, any $n$-dimensional Lie algebra $\frak{g}=\left(\mathbb{C}^{n},\mu\right)$ is identified with its law $\mu$ in $\mathcal{L}^{n}$. Moreover, we convene that nonwritten brackets are either zero or obtained by antisymmetry.

\section{Lie algebra models}

Let $\frak{L}^{n}$ be the variety of Lie algebra laws over $\mathbb{C}^{n}$.
As known, the linear group $GL\left(  n,\mathbb{C}\right)  $ acts on
$\mathbb{C}^{n}$ by changes of basis. Recall that a Lie algebra $\frak{g}%
=\left(  \mathbb{C}^{n},\lambda\right)  $ is said to contract to the law $\mu$
if $\mu\in\overline{O\left(  \lambda\right)  }$, where $\overline{O\left(
\lambda\right)  }$ is the closure of the orbit of $\lambda$ by the group
action ( in either the metric or Zariski topology [1]). It is easy to verify that this definition is quite the same as the more classical one: let $\left(\mathbb{C}^{n},\mu_{0}\right)$ be a point in $\mathcal{L}^{n}$
\cite{LeN} and $\left\{  f_{t}\right\}  $ a sequence of isomorphisms in
$\mathbb{C}^{n}$. Clearly the element $\left(  \mathbb{C}^{n},\mu_{t}\right)
$ with
\[
\mu_{t}=f^{-1}\circ\mu_{0}\left(  f_{t}\times f_{t}\right)
\]
belongs to the orbit $\mathcal{O}\left(  \mu_{0}\right)  $. If the limit
$\lim_{t\rightarrow\infty}\mu_{t}$ exists, it defines a law of $\mathcal{L}%
^{n}$ called a contraction of $\mu_{0}$. \newline As told before, for
analyzing Lie algebra models the concept of perturbation is more convenient.
This notion is better adapted for the study of certain topological properties
concerning the orbits of laws in the variety $\mathcal{L}^{n}$, and is placed
within the internal set theory ( I.S.T. \cite{Ne}).

\begin{definition}
A perturbation $\mu$ of a Lie algebra law $\mu_{0}$ in $\mathcal{L}^{n}$ is a
law satisfying%
\[
\mu\left(  X,Y\right)  \simeq\mu_{0}\left(  X,Y\right)  ,\;\;X,Y\in
\mathbb{C}^{n}\;\rm{standard}%
\]
i.e., the structure constants of $\mu$ are infinitely close to those of
$\mu_{0}$ ( the basis being fixed).
\end{definition}

Perturbations of a law can be decomposed according to the following rule \cite{G4}:%
\[
\mu=\mu_{0}+\varepsilon_{1}\varphi_{1}+...+\varepsilon_{1}..\varepsilon
_{k}\varphi_{k}%
\]
where $\varphi_{1},..,\varphi_{k}$ are (standard) $2$-cochains independent of
$\varepsilon_{1},..,\varepsilon_{k}\simeq0$. The integer $k$ is called the
length of the perturbation. If $k=1$, they correspond to the infinitesimal
deformations.\newline Let now $\frak{g}%
_{0}=\left(  \mathbb{C}^{n}.\mu_{0}\right)  $ be a Lie algebra and $\left(
P\right)  $ be an internal property concerning the Lie algebra laws.

\begin{definition}
The algebra $\frak{g}_{0}$ is called a semimodel relative to property $\left(
P\right)  $ if any Lie algebra law $\mu$ satisfying $\left(  P\right)  $
contracts on $\mu_{0}$.\newline The semimodel is called a model if, in
addition, any perturbation of $\mu_{0}$ satisfies $\left(  P\right)  $.
\end{definition}

\begin{remark}
The neighbourhoods of model algebras are entirely characterized by property
$\left(  P\right)  $ ( [2],
[3]). As examples of properties whose analysis lead to interesting questions
we could enumerate those related to differential invariants of forms and
systems on a Lie algebra.
\end{remark}

\begin{example}
Let $\frak{g}$ be a $\left(  2n+1\right)  $-dimensional Lie algebra and
$\omega$ a linear form. Then $\omega$ is called a contact form if
\ $\omega\wedge\left(  d\omega\right)  ^{n}\neq0$, where $d\omega\left(
X,Y\right)  =-\omega\left(  \left[  X,Y\right]  \right)  $ corresponds to the
contragradient representation of $\frak{g}$. Let $\left(  P\right)  $ be the
property ''there exists a linear contact form on $\frak{g}$''. It is not
difficult to see that the Heisenberg Lie algebra $\frak{h}_{n}$ given by
\[
\left\{  d\omega_{2n+1}=\sum_{k=0}^{n-1}\omega_{2k+1}\wedge\omega
_{2k+2}\right.
\]
is the unique model relative to this property \cite{G2}.
\end{example}

Another example illustrating the interest of models for physical applications is the following:

\begin{example}
Let $\frak{g}=\frak{so}(3)$. It can be easily seen, without using the fact that this algebra is rigid for being simple, that this algebra is the unique model for the property "any linear form on $\frak{g}$ is a contact form". This implies in particular that this property is valid only in dimension three \cite{G1}.
\end{example}

It could also happen that a model is not unique, even that a model for a
property $\left(  P\right)  $ does not exist ( as happens for the r-contact
systems [5]). Therefore the definition of (semi-)model has to be generalized.

\begin{definition}
A family $F$ of Lie algebras satisfying an internal property $\left(
P\right)  $ is called a multiple semimodel relative to $\left(  P\right)  $ if
any Lie algebra satisfying it contracts to an algebra of the family.\newline
If in addition any perturbation of elements in $F$ also satisfies $\left(
P\right)  $, the family is called a multiple model.
\end{definition}

In order to get a minimal family of models, we also define irreducible
multiple models:

\begin{definition}
Let $F=\left\{  \frak{g}_{i}=\left(  \mathbb{C}^{n},\mu_{i}\right)  \right\}
_{i\in I}$ be a multiple model relative to property $\left(  P\right)  $. Then
$F$ is called irreducible if for any $i,j\in I$ we have $\frak{\mu}_{i}%
\notin\overline{O\left(  \mu_{j}\right)  }$.
\end{definition}

The property we consider here is ''there exists a linear form $\omega$ whose
differential $d\omega$ is symplectic''. Recall that a linear form $\omega$ on
a $2n$-dimensional Lie algebra $\frak{g}$ is called symplectic if it is closed
and $\omega^{n}\neq0$. Considering the differentials $d\omega$ instead of
$\omega$ we obtain the frobeniusian Lie algebras:

\begin{definition}
A $2n$-dimensional Lie algebra $\frak{g}$ is called frobeniusian if there
exists $\omega\in\frak{g}^{\ast}$ such that $\left(  d\omega\right)  ^{n}%
\neq0$.
\end{definition}

Semimodels of frobeniusian Lie algebras were studied in [6] and [7]. Here
models exist, but uniqueness is lost. Let $\left(  \varphi\right)  :=\left(
\varphi_{1},..,\varphi_{p-1}\right)  \in\mathbb{C}^{p-1}$.

\begin{theorem}
\cite{G1} Let $\{F_{\varphi}\;|\;\left(  \varphi\right)  \in\mathbb{C}^{p-1}\}$
be the family on $\left(  p-1\right)  $-parameters of $2p$-dimensional Lie
algebras given by
\[
\left\{
\begin{array}
[c]{l}%
d\omega_{1}=\omega_{1}\wedge\omega_{2}+\sum_{k=1}^{p-1}\omega_{2k+1}%
\wedge\omega_{2k+2}\\
d\omega_{2}=0\\
d\omega_{2k+1}=\varphi_{k}\omega_{2}\wedge\omega_{2k+1},\;1\leq k\leq p-1\\
d\omega_{2k+2}=-\left(  1+\varphi_{k}\right)  \omega_{2}\wedge\omega
_{2k+2},\;1\leq k\leq p-1
\end{array}
\right.
\]
where $\left\{  \omega_{1},..,\omega_{2p}\right\}  $ is a basis of $\left(
\mathbb{C}^{2p}\right)  ^{\ast}$. The family $\{F_{\varphi}\}$ is an
irreducible multiple model for the property ''there exists a linear form whose
differential is symplectic''.
\end{theorem}

It can be easily seen that the algebras $F_{\varphi}$ admit the following graduation: if $\left\{X_{1},..,X_{2p}\right\}$ is a dual basis to $\left\{\omega_{1},..,\omega_{2p}\right\}$, then $F_{\varphi}= \left(F_{\varphi}\right)_{0}\oplus \left(F_{\varphi}\right)_{1}\oplus \left(F_{\varphi}\right)_{2}$, where $\left(F_{\varphi}\right)_{0}=\mathbb{C}{X_{2}}, \left(F_{\varphi}\right)_{1}=\sum_{k=3}^{2p}\mathbb{C}{X_{k}}$ and $\left(F_{\varphi}\right)_{2}=\mathbb{C}{X_{1}}$. This decompostion will be of importance for cohomological computations. 

\bigskip

Let $\Omega_{1}$ denote the union of the following subsets of $\mathbb{C}%
^{p-1}$
\begin{eqnarray*}
&  \left\{  1+\varphi_{i}+\varphi_{j}=0\right\}  _{1\leq i,j\leq p-1}\\
&  \left\{  2+\varphi_{i}+\varphi_{j}=0\right\}  _{1\leq i,j\leq p-1}\\
&  \left\{  \varphi_{i}-\varphi_{j}=0\right\}  _{1\leq i\neq j\leq p-1}\\
&  \left\{  \varphi_{i}=0\right\}  _{1\leq i\leq p-1}\\
&  \left\{  \varphi_{i}+1=0\right\}  _{1\leq i\leq p-1}\\
&  \left\{  2\varphi_{i}+1=0\right\}  _{1\leq i\leq p-1}%
\end{eqnarray*}

\begin{lemma}
If $\left(  \varphi\right)  \notin\Omega_{1}$, then $\dim Der\left(
F_{\varphi}\right)  =3p-1$.
\end{lemma}

\begin{proof}
Let $f\in Der\left(  F_{\varphi}\right)  $ and let $\left\{  X_{1}%
,..,X_{2p}\right\}  $ be a dual basis to $\left\{  \omega_{1},..,\omega
_{2p}\right\}  $. Let $f\left(  X_{i}\right)  =f_{i}^{j}X_{j}\;\left(  1\leq
i\leq2p\right)  $. For $1\leq k\leq p-1$ the conditions $\left[  f\left(
X_{2}\right)  ,X_{2k+1}\right]  +\left[  X_{2},f\left(  X_{2k+1}\right)
\right]  =\varphi_{k}f\left(  X_{2k+1}\right)  $\ imply $f_{2k+1}^{j}=0$ for
$j\neq1,2k+1$, as $\left(  \varphi\right)  \notin\Omega$, and $f_{2}%
^{2k+2}=-\left(  1+\varphi_{k}\right)  f_{2k+1}^{1}$. Similarly we obtain
$f_{2k+2}^{j}=0$ for $j\neq1,2k+2$ and $f_{2}^{2k+1}=-\varphi_{k}f_{2k+1}^{1}%
$. As $X_{2}$ does not belong to the derived subalgebra $D^{1}F_{\varphi}$, we have $f_{2}^{2}=0$. \newline
Therefore, any nontrivial derivation $f$ modulus the image $ad\left(
F_{\varphi}\right)  $ is given by
\[
f\left(X_{i}\right)  \;\left(\it{mod}\;ad\left(  F_{\varphi
}\right)  \right)  =\left\{
\begin{tabular}
[c]{ll}%
$0$ & $i=1,2$\\
$\lambda_{i}X_{i}$ & $i$ odd and $\geq3$\\
$-\lambda_{i}X_{i}$ & $i$ even and $\geq4$%
\end{tabular}
\right.
\]
where $\lambda_{i}\in\mathbb{C}$. This is a linear combination of the
derivations $f_{i}\;\left(  1\leq i\leq p-1\right)  $ defined by
\[
\left\{
\begin{array}
[c]{c}%
f_{i}\left(  X_{2i+1}\right)  =X_{2i+1}\\
f_{i}\left(  X_{2i+2}\right)  =-X_{2i+2}%
\end{array}
\right.
\]
This shows that the $\left\{  f_{i}\right\}  _{1\leq i\leq p-1}$ form a basis
of $H^{1}\left(  F_{\varphi},F_{\varphi}\right)  $, thus dim $H^{1}\left(
F_{\varphi},F_{\varphi}\right)  =p-1$. The algebra $Der\left(  F_{\varphi
}\right)  $ is easily seen to be isomorphic to the semidirect product
$adF_{\varphi}\oplus\sum_{i=1}^{p-1}\mathbb{C}\left\{  f_{i}\right\}  $ given
by
\[
\left\{  \left[  ad\left(  X_{2i+1+t}\right)  ,f_{i}\right]  =-\left(
-1\right)  ^{t}ad\left(  X_{2i+1+t}\right)  ,\;1\leq i\leq p-1,\;t\in\left\{
0,1\right\}  \right.
\]
Clearly $\dim\,Der\left(  F_{\varphi}\right)  =3p-1$.
\end{proof}

\begin{corollary}
If $p\geq2$, then $dim H^{1}\left(  F_{\varphi},F_{\varphi}\right)  =p-1$.
\end{corollary}

Now let $\frak{g}$ be a frobeniusian Lie algebra not belonging to the family
$\{F_{\varphi}\;|\;\varphi\in\mathbb{C}^{p-1}\}$. Then $\frak{g}$ contracts to
some element $F_{\varphi_{0}}$ in view of theorem 1. We denote it by
$\frak{g}\rightarrowtail F_{\varphi_{0}}$. By the general properties of
contractions \cite{Bu}, $\frak{g}$ satisfies the following conditions

\begin{enumerate}
\item $\dim Der\left(  \frak{g}\right)  <\dim Der\left(  F_{\varphi_{0}}\right)  $

\item $\dim\left[  \frak{g},\frak{g}\right]  \geq\dim\left[  F_{\varphi_{0}
},F_{\varphi_{0}}\right]  $

\item $\dim Z\left(  \frak{g}\right)  \leq\dim Z\left(  F_{\varphi_{0}}\right)  $
\end{enumerate}

where $Der\left(  \circ\right)  $ denotes the algebra of derivations, $\left[
\circ,\circ\right]  $ the commutator algebra and $Z\left(  \circ\right)  $ the
center. As the algebras $F_{\varphi}$ are centerless, we have $Z\left(
\frak{g}\right)  =\left\{  0\right\}  $ for a frobeniusian Lie algebra.

\begin{proposition}
Let $p\geq2$ and $\frak{g}$ be a frobeniusian Lie algebra
$\frak{g\rightarrowtail} F(\varphi)$ such that $\left(  \varphi\right)
\notin\Omega_{1}$. Then $\dim H^{1}\left(  \frak{g},\frak{g}\right)  \leq p-2$.
\end{proposition}

\begin{proof}
As $Z\left(  \frak{g}\right)  =\left\{  0\right\}  ,\;$the adjoint
representation is faithful and we obtain
\[
\dim H^{1}\left(  \frak{g},\frak{g}\right)  =\dim Der\left(  \frak{g}\right)
-2p<\dim Der\left(  F_{\varphi}\right)  -2p
\]
The result follows from the previous lemma.
\end{proof}

One of the reasons for isolating the parameters of $\Omega$ is the similar structure of its derivations, which leads to isomorphisms of the corresponding derivation algebras. In particular, we will obtain complete Lie algebras. Recall that a Lie algebra is called complete if $Z\left(  \frak{g}\right)
=H^{1}\left(  \frak{g},\frak{g}\right)  =\left\{  0\right\}  $.

\begin{theorem}
Let $p\geq2$. If $\left(  \varphi\right)  \allowbreak\notin\Omega_{1}$, the Lie
algebra $Der\left(  F_{\varphi}\right)$ is a $3$-step solvable complete Lie algebra. 
Moreover, for any $\left(  \varphi_{0}\right)  \neq\left(  \varphi_{0}^{\prime}\right)
$ we have
\[
Der\left(  F_{\varphi_{0}}\right)  \simeq Der\left(  F_{\varphi
_{0}^{\prime}}\right)
\]
\end{theorem}

\begin{proof}
The structural equations of $Der\left(  F_{\varphi}\right)  $ for any $\left(
\varphi\right)  \notin\Omega_{1}$ are given by%
\[
\left\{
\begin{array}
[c]{l}%
d\omega_{1}=\omega_{1}\wedge\omega_{2}+\sum_{k=1}^{p-1}\omega_{2k+1}%
\wedge\omega_{2k+2}\\
d\omega_{2}=0\\
d\omega_{2k+1}=\varphi_{k}\omega_{2}\wedge\omega_{2k+1}+\omega_{2p+k}%
\wedge\omega_{2k+1},\;1\leq k\leq p-1\\
d\omega_{2k+2}=-\left(  1+\varphi_{k}\right)  \omega_{2}\wedge\omega
_{2k+2}-\omega_{2p+k}\wedge\omega_{2k+2},\;1\leq k\leq p-1\\
d\omega_{2p+k}=0,\;1\leq k\leq p-1
\end{array}
\right.
\]
over the basis $\left\{  \omega_{1},..,\omega_{2p},\omega_{2p+1}%
,..,\omega_{3p-1}\right\}  $. The proof of completeness is routine, as well as
the fact that the algebra is 3-step solvable. For the last assertion, let
$\left(  \varphi_{0}\right)  \neq\left(  \varphi_{0}^{\prime}\right)  $ and
consider the change of basis given by%
\[
\left\{
\begin{array}
[c]{l}%
X_{2}^{\prime}=X_{2}+\displaystyle\sum_{i=1}^{p-1}\left(  \varphi_{i}^{\prime}-\varphi
_{i}\right)  X_{2p+i}\\
X_{i}^{\prime}=X_{i},\;i\neq2
\end{array}
\right.
\]
\end{proof}

\begin{remark}
Observe in particular that $Der\left(  F_{\varphi}\right)  $ cannot arise
as a contraction of a Lie algebra, since it is complete \cite{Car}.
 Moreover, it is easily seen that this Lie algebra is rigid \cite{AC}, as it is isomorphic to the
semidirect product $\frak{h}_{p-1}\oplus\frak{t}$ of the $\left(  2p-1\right)
$-dimensional Heisenberg Lie algebra $\frak{h}_{p-1}$ with a maximal torus of derivations.
\end{remark}

\section{Chevalley cohomology}

In this section we determine the cohomology groups $H^{2}\left(  F_{\varphi
},F_{\varphi}\right)  $ whenever $\varphi$ lies in some dense open subset of
$\mathbb{C}^{p-1}$. Specifically, let $\Omega_{2}$ be the union of following
hyperplanes in $\mathbb{C}^{p-1}:$%
\begin{eqnarray*}
&  \left\{  1+\varphi_{i}-\varphi_{j}=0\right\}  \\
&  \left\{  \varphi_{i}+\varphi_{j}=0\right\}  \\
&  \left\{  2+\varphi_{i}=0\right\}  \\
&  \left\{  1-\varphi_{i}=0\right\}  \\
&  \left\{  1+2\varphi_{i}-\varphi_{j}=0\right\}  \\
&  \left\{  1+2\varphi_{i}+\varphi_{j}=0\right\}  \\
&  \left\{  2\varphi_{j}-\varphi_{j}=0\right\}  \\
&  \left\{  2+2\varphi_{i}+\varphi_{j}=0\right\}
\end{eqnarray*}
and $\Omega=\Omega_{1}\cup\Omega_{2}$. Clearly this is a finite union of
hyperplanes, and therefore its complementary in $\mathbb{C}^{p-1}$ is dense.
In what follows we always suppose that $\left(  \varphi\right)  \notin\Omega$.
In the previous section we have seen that the algebras $F_{\varphi}$ are graded. Therefore,
the cohomology spaces inherit a graduation ( see \cite{Ko}, as well as for the notation used), and we have
$H^{2}\left(  F_{\varphi},F_{\varphi}\right)  =F_{-k}H^{2}\left(  F_{\varphi
},F_{\varphi}\right)  $ for some positive integer $k$, where 
$F_{-k}H^{2}\left(  F_{\varphi},F_{\varphi}\right)  =\sum_{j\geq -k}^{1}H_{k}%
^{2}\left(  F_{\varphi},F_{\varphi}\right)$. For brevity in the notation, we will denote the coboundary operator $\delta_{\mu_{\varphi},p}:C^{p}\left(F_{\varphi},F_{\varphi}\right)\rightarrow C^{p+1}\left(F_{\varphi},F_{\varphi}\right)$ simply by $\delta_{\mu_{\varphi}}$

\begin{lemma}
For $p\geq2$ we have $H^{2}\left(  F_{\varphi},F_{\varphi}\right)
=F_{-2}H^{2}\left(  F_{\varphi},F_{\varphi}\right)$.
\end{lemma}

\begin{proof}
Let $\psi\in Z^{2}\left(  F_{\varphi},F_{\varphi}\right)  $ with $\psi\left(
X_{i},X_{j}\right)  =\sum_{k}a_{ij}^{k}X_{k}$ $\left(  1\leq i<j\leq2p\right)
$. Considering the coboundary operator $\delta_{\mu_{\varphi}}$ for the
triples $\left\{  X_{1},X_{2k+1},X_{2k+2}\right\}  _{1\leq k\leq p-1}$ we
obtain%
\begin{eqnarray*}
\delta_{\mu_{\varphi}}\left(  \psi\right)  \left(  X_{1},X_{2k+1}%
,X_{2k+2}\right)    & =-a_{1,2k+1}^{2}\left(  1+\varphi_{k}\right)
X_{2k+2}-a_{1,2k+2}^{2}\varphi_{k}X_{2k+1}+\\
+\sum_{i\neq2k+1,2k+2}\alpha_{i}X_{i},\;  & (\alpha_{i}\in\mathbb{C})%
\end{eqnarray*}
Since $\left(  \varphi\right)  \notin\Omega$ we have $\varphi_{k}\left(
1+\varphi_{k}\right)  \neq0$ for any $k$, and thus $a_{1,2k+1}^{2}%
=a_{1,2k+2}^{2}=0\;\left(  1\leq k\leq p-1\right)  $. This shows that
$Z_{-3}^{2}\left(  F_{\varphi},F_{\varphi}\right)  =0$. Now $Z_{k}^{2}\left(
F_{\varphi},F_{\varphi}\right)  =0$ for $k\geq2$, as follows immediately from
the graduation.
\end{proof}

Therefore, the determination of the cohomology spaces reduces to the computation of the distinct subspaces $H^{2}_{k}\left(F_{\varphi},F_{\varphi}\right)$. We begin determining bases for the cocycle spaces $Z^{2}_{k}\left(F_{\varphi},F_{\varphi}\right)$.

\begin{lemma}
$\dim$ $Z_{-2}^{2}\left(  F_{\varphi},F_{\varphi}\right)  =1$.
\end{lemma}

\begin{proof}
Let $\psi\in Z_{-2}^{2}\left(  F_{\varphi},F_{\varphi}\right)  $. Then its
generic form is
\begin{eqnarray*}
\psi\left(  X_{1},X_{2}\right)   &  =a_{12}^{2}X_{2}\\
\psi\left(  X_{i},X_{j}\right)   &  =a_{ij}^{2}X_{2},\;3\leq i,j\leq2p\\
\psi\left(  X_{1},X_{j}\right)   &  =\sum_{t=3}^{2p}a_{1,j}^{t}X_{t},\;j\geq3
\end{eqnarray*}
For $p\geq4$ we can always find $3\leq l\leq2p$ such that $[X_{i}%
,X_{j}]=[X_{i},X_{k}]=[X_{j},X_{k}]=0$. Then
\[
\delta_{\mu_{\varphi}}\left(  \psi\right)  \left(  X_{i},X_{j},X_{l}\right)
=a_{ij}^{2}[X_{2},X_{l}]+\sum_{t\neq l}\alpha_{t}X_{t},\;\alpha_{t}%
\in\mathbb{C}%
\]
where
\[
\left[  X_{2},X_{l}\right]  =\left\{
\begin{array}
[c]{l}%
\varphi_{k}X_{2k+1}\;\rm{if }l=2k+1\\
-\left(  1+\varphi_{k}\right)  X_{2k+2}\;\rm{if\ }l=2k+2
\end{array}
\right.
\]
This implies $a_{ij}^{2}=0$ whenever $[X_{i},X_{j}]=0$. For $p=2,3$ the nullity follows from considering the triples $\left\{  X_{2},X_{i},X_{j}\right\}
$.\newline Applying $\delta_{\mu_{\varphi}}$ to the triples $\left\{
X_{2},X_{2k+1},X_{2k+2}\right\}  $ we obtain
\[
a_{12}^{2}=a_{2k+1,2k+2}^{2},\;1\leq k\leq p-1
\]
Finally, taking $\left\{  X_{1},X_{2},X_{j}\right\}  _{j\geq3}$ we get%
\begin{eqnarray*}
a_{1,j}^{k} &  =0,\rm{\ }k\neq j\\
a_{1,j}^{j}+\beta_{j}a_{12}^{2} &  =0
\end{eqnarray*}
Thus $Z_{-2}^{2}\left(  F_{\varphi},F_{\varphi}\right)  $ is generated by the
cocycle $\psi_{12}^{2}$ defined by

\[
\left\{
\begin{tabular}
[c]{ll}%
$\left(  X_{1},X_{2}\right)  \mapsto X_{2}$ & \\
$\left(  X_{1},X_{2k+1}\right)  \mapsto -\varphi_{k}X_{2k+1},\;1\leq k\leq p-1$ & \\
$\left(  X_{1},X_{2k+2}\right)  \mapsto (1+\varphi_{k})X_{2k+2},\;1\leq k\leq p-1$ & \\
$\left(  X_{2k+1},X_{2k+2}\right)  \mapsto X_{2},\;1\leq k\leq p-1  $ & %
\end{tabular}
\right. 
\]
\end{proof}

\begin{lemma}
$\dim\;Z_{-1}^{2}\left(  F_{\varphi},F_{\varphi}\right)  =4\left(  p-1\right)  .$
\end{lemma}

\begin{proof}
A cocycle $\psi\in Z_{-1}^{2}\left(  F_{\varphi},F_{\varphi}\right)  $ has the
generic form%
\begin{eqnarray*}
\psi\left(  X_{2},X_{j}\right)   &  =a_{2j}^{2}X_{2},\;j\geq3\\
\psi\left(  X_{1},X_{2}\right)   &  =\sum_{j=3}^{2p}a_{12}^{j}X_{j}\\
\psi\left(  X_{i},X_{j}\right)   &  =\sum_{t=3}^{2p}a_{ij}^{t}X_{t},\;3\leq
i,j\leq2p\\
\psi\left(  X_{1},X_{j}\right)   &  =a_{1j}^{1}X_{1},\;3\leq j\leq2p
\end{eqnarray*}
and satisfies the following equations, for any $1\leq k\leq p-1:$%
\begin{eqnarray*}
a_{12}^{2k+2-r}+\left(  -1\right)  ^{r}a_{2,2k+1+r}^{2}+\left(  r+\varphi
_{k}\right)  a_{1,2k+1+r}^{1}  & =0\\
a_{12}^{2k+2-r}-\left(  1+\varphi_{k}-r\right)  a_{2,2k+1+r}^{2}+\left(
-1\right)  ^{r}\left(  r+\varphi_{k}\right)  a_{2k+1,2k+2}^{2k+2-r}  & =0\\
a_{12}^{2k+2-r}+\left(  -1\right)  ^{r}(r+\varphi_{j})a_{2j+1,2j+2}^{2k+2-r}
& =0\\
\varphi_{j}a_{2,2k+1+r}^{2}+\left(  -1\right)  ^{r}(r+\varphi_{k}%
)a_{2k+1+r,2j+1}^{2j+1}  & =0\\
\left(  1+\varphi_{j}\right)  a_{2,2k+1+r}^{2}-\left(  -1\right)
^{r}(r+\varphi_{k})a_{2k+1+r,2j+2}^{2j+2}  & =0
\end{eqnarray*}
where $1\leq j\neq k\leq p-1$ and  $r\in\left\{  0,1\right\}  $. From the
structure of this system it is not difficult to extract a basis. One is given
by the cocycles $\psi_{2,2k+1}^{2},\psi_{2,2k+2}^{2},\psi_{12}^{2k+1}$%
,\newline $\psi_{12}^{2k+2}$ $\left(  1\leq k\leq p-1\right)  $ defined by :

\begin{enumerate}
\item $\psi_{2,2k+1}^{2}$:
\[
\left\{  
\begin{tabular}
[c]{ll}%
$\left(  X_{1},X_{2k+1}\right)  \mapsto X_{1};\;\;\left(  X_{2},X_{2k+1}%
\right)  \mapsto-\varphi_{k}X_{2}$ & \\
$\left(  X_{2k+1},X_{2k+2}\right)  \mapsto-\left(  1+\varphi_{k}\right)
X_{2k+2}$ & \\
$\left(  X_{2k+1},X_{2k^{\prime}+1}\right)  \mapsto\varphi_{k^{\prime}%
}X_{2k^{\prime}+1}\;$ & $\left(  k\neq k^{\prime}\right)  $\\
$\left(  X_{2k+1},X_{2k^{\prime}+2}\right)  \mapsto-\left(  1+\varphi
_{k^{\prime}}\right)  X_{2k^{\prime}+2}$ & $\left(  k\neq k^{\prime}\right)  $%
\end{tabular}
\right.  
\]

\item $\psi_{2,2k+2}^{2}$:
\[
\left\{
\begin{tabular}
[c]{ll}%
$\left(  X_{1},X_{2k+2}\right)  \mapsto X_{1};\;\;\left(  X_{2},X_{2k+2}%
\right)  \mapsto(1+\varphi_{k})X_{2}$ & \\
$\left(  X_{2k+1},X_{2k+2}\right)  \mapsto-\varphi_{k}X_{2k+1}$ & \\
$\left(  X_{2k+2},X_{2k^{\prime}+1}\right)  \mapsto\varphi_{k^{\prime}%
}X_{2k^{\prime}+1}\;$ & $\left(  k\neq k^{\prime}\right)  $\\
$\left(  X_{2k+2},X_{2k^{\prime}+2}\right)  \mapsto-\left(  1+\varphi
_{k^{\prime}}\right)  X_{2k^{\prime}+2}$ & $\left(  k\neq k^{\prime}\right)  $%
\end{tabular}
\right. 
\]

\item $\psi_{12}^{2k+1}$:
\[
\left\{
\begin{tabular}
[c]{ll}%
$\left(  X_{1},X_{2}\right)  \mapsto\left(  1+\varphi_{k}\right)  X_{2k+1}$ &
\\
$\left(  X_{1},X_{2k+2}\right)  \mapsto-X_{1}$ & \\
$\left(  X_{2j+1},X_{2j+2}\right)  \mapsto X_{2k+1}$ & $;1\leq j\leq p-1$%
\end{tabular}
\right.  
\]

\item $\psi_{12}^{2k+2}$:
\[
\left\{
\begin{tabular}
[c]{ll}%
$\left(  X_{1},X_{2}\right)  \mapsto\varphi_{k}X_{2k+2}$ & \\
$\left(  X_{1},X_{2k+1}\right)  \mapsto-X_{1}$ & \\
$\left(  X_{2j+1},X_{2j+2}\right)  \mapsto-X_{2k+2}$ & $;1\leq j\leq p-1$%
\end{tabular}
\right.  
\]
\end{enumerate}
Therefore $\dim\;Z_{-1}^{2}\left(  F_{\varphi},F_{\varphi}\right)  =4\left(
p-1\right)$.  
\end{proof}

\begin{lemma}
$\dim\;Z_{0}^{2}\left(  F_{\varphi},F_{\varphi}\right)  =4p^{2}-8p+5.$
\end{lemma}

\begin{proof}
Any cocycle $\psi$ belonging to this space has the generic form
\begin{eqnarray*}
\psi\left(  X_{2},X_{2k+1+r}\right)   &  =\sum_{j=3}^{2p}a_{2,2k+1+r}^{j}%
X_{j},\;r=0,1\\
\psi\left(  X_{1},X_{2}\right)   &  =a_{12}^{1}X_{1}\\
\psi\left(  X_{i},X_{j}\right)   &  =a_{ij}^{1}X_{1},\;3\leq i,j\leq2p
\end{eqnarray*}
Now, by the coboundary operator $\delta_{\mu_{\varphi}}$ the following
equations must be satisfied:%
\begin{eqnarray*}
a_{12}^{1}+a_{2,2k+1}^{2k+1}+a_{2,2k+2}^{2k+2} &  =0,\;1\leq k\leq p-1\\
a_{2,2k+1}^{2t+1}+\left(  \varphi_{k}-\varphi_{t}\right)  a_{2k+1,2t+2}%
^{1}+a_{2,2t+2}^{2k+2} &  =0,\;1\leq k\neq t\leq p-1\\
a_{2,2t+1}^{2k+2}+\left(  1+\varphi_{t}+\varphi_{k}\right)  a_{2k+1,2t+1}%
^{1}-a_{2,2k+1}^{2t+2} &  =0,\;1\leq k\neq t\leq p-1\\
a_{2,2k+2}^{2t+1}-\left(  1+\varphi_{t}+\varphi_{k}\right)  a_{2k+2,2t+2}%
^{1}-a_{2,2t+2}^{2k+1} &  =0,\;1\leq k\neq t\leq p-1
\end{eqnarray*}
From the equations we deduce the following linearly independent cocycles:

\begin{enumerate}
\item $\psi_{12}^{1}$:
\[
\left\{
\begin{tabular}
[c]{ll}%
$\left(  X_{1},X_{2}\right)  \mapsto X_{1}$ & \\
$\left(  X_{2},X_{2k+2}\right)  \mapsto-X_{2k+2}$ & $,1\leq k\leq p-1$%
\end{tabular}
\right.  
\]

\item $\psi_{2,2k+1}^{2k+1}\;\left(  1\leq k\leq p-1\right)  :$%
\[
\left\{
\begin{tabular}
[c]{ll}%
$\left(  X_{2},X_{2k+1}\right)  \mapsto X_{2k+1},$ & $1\leq k\leq p-1$\\
$\left(  X_{2},X_{2k+2}\right)  \mapsto-X_{2k+2},$ & $1\leq k\leq p-1$%
\end{tabular}
\right.
\]

\item $\psi_{2,2k+1}^{2t+2}\;\left(  1\leq k\leq p-1\right)  :$%
\[
\left\{
\begin{tabular}
[c]{ll}%
$\left(  X_{2},X_{2k+1}\right)  \mapsto-\left(  1+\varphi_{t}+\varphi
_{k}\right)  X_{2t+2},$ & $1\leq k\neq t\leq p-1$\\
$\left(  X_{2k+1},X_{2t+1}\right)  \mapsto-X_{1}$ & $1\leq k\neq t\leq p-1$%
\end{tabular}
\right.
\]

\item $\psi_{2,2k+2}^{2t+1}\;\left(  1\leq k\leq p-1\right)  :$%
\[
\left\{
\begin{tabular}
[c]{ll}%
$\left(  X_{2},X_{2k+2}\right)  \mapsto\left(  1+\varphi_{t}+\varphi
_{k}\right)  X_{2t+1},$ & $1\leq k\neq t\leq p-1$\\
$\left(  X_{2k+2},X_{2t+2}\right)  \mapsto X_{1}$ & $1\leq k\neq t\leq p-1$%
\end{tabular}
\right.
\]

\item $\psi_{2,2k+2}^{2t+2}\;\left(  1\leq k\leq p-1\right)  :$%
\[
\left\{
\begin{tabular}
[c]{ll}%
$\left(  X_{2},X_{2k+2}\right)  \mapsto(\varphi_{k}-\varphi_{t})X_{2t+2},$ &
$1\leq k\neq t\leq p-1$\\
$\left(  X_{2k+2},X_{2t+1}\right)  \mapsto-X_{1}$ & $1\leq k\neq t\leq p-1$%
\end{tabular}
\right.
\]
To these we have to add those cocycles for which the operator
$\delta_{\mu_{\varphi}}$ gives no conditions:

\item $\psi_{2k+1,2k+2}^{1}\;\left(  1\leq k\leq p-1\right)  :$
\[
\left\{  \left(
X_{2k+1},X_{2k+2}\right)  \mapsto X_{1}\right.  
\]

\item $\psi_{2,2k+1}^{2k+2}\;\left(  1\leq k\leq p-1\right) :$
\[
\left\{
\left(  X_{2},X_{2k+1}\right)  \mapsto X_{2k+2}\right.  
\]

\item $\psi_{2,2k+2}^{2k+1}\;\left(  1\leq k\leq p-1\right)  :$
\[
\left\{
\left(  X_{2},X_{2k+2}\right)  \mapsto X_{2k+1}\right.  
\]
\end{enumerate}

Adding up, we obtain
\[
\dim\;Z_{0}^{2}\left(  F_{\varphi},F_{\varphi}\right)  =1+4\left(  p-1\right)
+4\left(  p-1\right)  \left(  p-2\right)  =4p^{2}-8p+5
\]
\end{proof}

\begin{lemma}
$\dim\;Z_{1}^{2}\left(  F_{\varphi},F_{\varphi}\right)  =2\left(  p-1\right)  .$
\end{lemma}

\begin{proof}
These cocycles have the generic form $\psi\left(  X_{2},X_{j}\right)
=a_{2,j}X_{1}$ for $3\leq j\leq2p.$ The coboundary operator gives no
conditions on the coefficients, thus we have the basis \newline $\psi
_{2,j}^{1}\;\left(  3\leq j\leq2p\right)  :$%
\[
\left\{  \left(  X_{2},X_{j}\right)  \mapsto X_{1}\right.
\]
of dimension $2\left(  p-1\right)  $.
\end{proof}

\begin{proposition}
If $\left(  \varphi\right)  \notin\Omega$ then $\dim\;H^{2}\left(  F_{\varphi
},F_{\varphi}\right)  =p-1$.
\end{proposition}

\begin{proof}
By the preceding lemmas $\dim\;Z^{2}\left(  F_{\varphi},F_{\varphi}\right)
=4p^{2}-2p$ . Now $\dim\;B^{2}\left(  F_{\varphi},F_{\varphi}\right)
=4p^{2}-\left(  3p-1\right)$, since $\dim Der\left(  F_{\varphi}\right)
=3p-1$ by lemma 1.
\end{proof}

\begin{proposition}
For $p\geq2$ the following holds:

\begin{enumerate}
\item $B_{-2}^{2}\left(  F_{\varphi},F_{\varphi}\right)  =Z_{-2}^{2}\left(
F_{\varphi},F_{\varphi}\right)  $

\item $B_{-1}^{2}\left(  F_{\varphi},F_{\varphi}\right)  =Z_{-1}^{2}\left(
F_{\varphi},F_{\varphi}\right)  $

\item $\dim Z_{0}^{2}\left(  F_{\varphi},F_{\varphi}\right)  -\dim B_{0}%
^{2}\left(  F_{\varphi},F_{\varphi}\right)  =p-1$

\item $Z_{1}^{2}\left(  F_{\varphi},F_{\varphi}\right)  =B_{1}^{2}\left(
F_{\varphi},F_{\varphi}\right)  $
\end{enumerate}
\end{proposition}

\begin{proof}
The proof of \textit{(i),(ii)} and \textit{(iv)} follows immediately by application of the
coboundary operator to an element $f\in GL\left(  2p,\mathbb{C}\right)  $. Let
$\theta\in B_{0}^{2}\left(  F_{\varphi},F_{\varphi}\right)  $ be the
coboundary given by
\[
\left\{
\begin{array}
[c]{l}%
\theta\left(  X_{1},X_{2}\right)  =X_{1}\\
\theta\left(  X_{2},X_{2k+1}\right)  =\varphi_{k}X_{2k+1},\;1\leq k\leq p-1\\
\theta\left(  X_{2},X_{2k+2}\right)  =-\left(  1+\varphi_{k}\right)
X_{2k+2},\;1\leq k\leq p-1
\end{array}
\right.
\]
It is a straightforward verification that
\[
\psi_{12}^{1}+\sum_{k=1}^{p-1}\psi_{2,2k+1}^{2k+1}=\theta
\]
Thus a basis of $H^{2}\left(  F_{\varphi},F_{\varphi}\right)  $ is given by
the classes $\left[  \psi_{2,2k+1}^{2k+1}\right]  $ for $1\leq k\leq p-1$.
\end{proof}

Recall that the mapping that associates the product $\psi\circ\psi$ to a $2$-cocycle $\psi$ induces a map $sq_{1}:H^{2}\left(  F_{\varphi},F_{\varphi
}\right)  \rightarrow H^{3}\left(  F_{\varphi},F_{\varphi}\right)$ called the Rim map \cite{Rim}. The image $sq_{1}\left(\psi\right)$ gives the first obstruction to the integability of the class of $\psi$.

\begin{corollary}
For any $p\geq2$ the Rim map $sq_{1}:H^{2}\left(  F_{\varphi},F_{\varphi
}\right)  \rightarrow H^{3}\left(  F_{\varphi},F_{\varphi}\right)  $ is
identically zero.
\end{corollary}

\begin{proof}
For the representatives $\psi_{2,2k+1}^{2k+1}\;\left(  1\leq k\leq p-1\right)
$ of the nontrivial cohomology classes we have $\psi_{2,2k+1}^{2k+1}\circ
\psi_{2,2k+1}^{2k+1}=0$ since $X_{2}\notin\psi_{2,2k+1}^{2k+1}\left(  F_{\varphi
},F_{\varphi}\right)  $. Therefore the mapping is everywhere zero.
\end{proof}

It can also be easily seen that the cocycles $\psi_{2,2k+1}^{2k+1}$ are linearly expandable, i.e., they define a linear deformation $F_{\varphi}+t\psi_{2,2k+1}^{2k+1}$. The importance of
this fact follows from the next

\begin{theorem}
If $\left(  \varphi\right)  \notin\Omega$, almost any infinitesimal
deformation of $F_{\varphi}$ lies in the set $\left\{  F_{\varphi}\;|\;\left(
\varphi\right)  \notin\Omega\right\}  $.
\end{theorem}

\begin{proof}
We have seen that any nontrivial cohomology class admits a linearly expandable
representative. For any $1\leq k\leq p-1$ we have
\[
F_{\varphi}+t\psi_{2,2k+1}^{2k+1}\simeq F\left(  \varphi_{1},..,\varphi
_{k}+t,..,\varphi_{p-1}\right)
\]
and for small $t$ we have $\left(  \varphi_{1},..,\varphi_{k}+t,..,\varphi
_{p-1}\right)  \notin\Omega$.
\end{proof}

\begin{corollary}
If $\left(  \varphi\right)  \notin\Omega$, any linear deformation $F_{\varphi
}+t\psi_{2,2k+1}^{2k+1}$ belongs to the family $\left\{  F_{\varphi
}\;|\;\left(  \varphi\right)  \in\mathbb{C}^{p-1}\right\}  $.
\end{corollary}

This corollary shows that frobeniusian model Lie algebras admitting deformations not belonging to the model are relatively scarce. 

\begin{corollary}
Let $\frak{g}$ be a non-model  frobeniusian Lie algebra. Then $\frak{g}$
contracts on some model Lie algebra $F_{\varphi}$ \ for $\left(
\varphi\right)  \in\Omega$. 
\end{corollary}

\begin{remark}
In particular, from this corollary we deduce that there exists only a finite
number of hyperplanes in $\mathbb{C}^{p-1}$ such that the algebras
$F_{\varphi}$ whose parameters  $\left(  \varphi\right)  $ belong to these
spaces allow deformations which lie outside the multiple model.
\end{remark}

\begin{corollary}
If $\frak{g}\notin$ $\left\{  F_{\varphi}\;|\;\left(  \varphi\right)
\notin\Omega\right\}  $ is frobeniusian, then $\dim\;H^{2}\left(
F_{\varphi_{0}},F_{\varphi_{0}}\right)  >p-1$ for any $F_{\varphi_{0}}$ such
that $F_{\varphi_{0}}\in\overline{\mathcal{O}\left(  \frak{g}\right)  }$. 
\end{corollary}

\begin{conclusion}
It has been shown that most of the Lie algebras admitting a linear form whose differential is symplectic already belong to the model, which constitutes a surprising result when compared with other properties such as the $r$-contact systems. The cohomology of the family $\left\{F_{\varphi}\right\}$ shows that the search for non-model frobeniusian Lie algebras reduces to study the parameters lying in $\Omega$. This implies that these algebras depend on less than $p-1$ parameters, due to the relations they satisfy. This translates in either more generators of the corresponding Lie algebra or eigenvalues of multiplicity $r\geq 2$ for the semisimple derivations of the algebra. Due to these particularities, their cohomology has to be computed for each case, since the structure depends heavily on the values the parameters take within $\Omega$. 
\end{conclusion}

\end{document}